\documentclass{article}
\usepackage{amsmath,amsthm,amssymb,authblk}

\newtheorem{theorem}{Theorem}[section]
\newtheorem{corollary}[theorem]{Corollary}

\newtheorem{lemma}[theorem]{Lemma}

\newtheorem{remark}[theorem]{Remark}
\numberwithin{equation}{section}
\def\pf{\noindent{\bf Proof.}~}

\begin{document}
\title{Permutation-like Matrix Groups with a Maximal Cycle
of Length Power of Two}
\author{Guodong Deng,\quad Yun Fan\\
\small School of Mathematics and Statistics\\
\small Central China Normal University, Wuhan 430079, China}
\date{}
\maketitle

\insert\footins{\footnotesize{\it Email address}:
yfan@mail.ccnu.edu.cn (Yun Fan).}

\begin{abstract}
If every element of a matrix group is similar to a permutation matrix,
then it is called a permutation-like matrix group. References \cite{C07}, \cite{DF}
and \cite{DF15} showed that, if a permutation-like matrix group
contains a maximal cycle such that the maximal cycle generates a normal subgroup
and the length of the maximal cycle equals to a prime, or a square of a prime,
or a power of an odd prime, then the permutation-like matrix group
is similar to a permutation matrix group.
In this paper, we prove that if a permutation-like matrix group
contains a maximal cycle such that the maximal cycle generates a normal subgroup
and the length of the maximal cycle equals to any power of $2$,
then it is similar to a permutation matrix group.


\end{abstract}

\section{Introduction}

By ${\rm GL}_{d}({\Bbb C})$ we denote the 
complex general linear group of dimension (degree)~$d$
consisting of invertible $d\times d$ complex matrices.
Any subgroup ${\cal G}$ of ${\rm GL}_{d}({\Bbb C})$ 
is said to be a {\em matrix group} of dimension $d$. 
If there is a $T\in{\rm GL}_d({\Bbb C})$
such that for every $A\in{\cal G}$ the conjugation 
$T^{-1}AT$ is a permutation matrix, 
then we say that ${\cal G}$ is a {\em permutation matrix group}. 
If for any $A\in{\cal G}$ there is a $T\in{\rm GL}_d({\Bbb C})$
such that $T^{-1}AT$ is a permutation matrix,
then ${\cal G}$ is called a {\em permutation-like matrix group}.
 
Cigler \cite{C05,C07} showed that a permutation-like matrix group
is not a permutation matrix group in general.
A $d\times d$ matrix is called a {\em maximal cycle} if
it is similar to a permutation matrix corresponding
to the cycle permutation of length~$d$.
Cigler conjectured that:

\medskip\hangindent15ex
\noindent{\bf Conjecture:}\,
{\it
 If a permutation-like matrix group contains a maximal cycle, then
it is a permutation matrix group.}

\medskip
Let ${\cal G}$ be a permutation-like matrix group 
containing a maximal cycle which generates a normal cyclic subgroup.
Cigler \cite{C05,C07} proved that if the dimension of ${\cal G}$ is a prime 
then ${\cal G}$ is a permutation matrix group.
In \cite{DF, DF15}, we further proved that 
the ${\cal G}$ is a permutation matrix group
if the dimension is a square of a prime, or 
any power of an odd prime.

In this paper we prove:

\begin{theorem}\label{main}
Let ${\cal G}$ be a permutation-like matrix group of dimension $2^n$
where $n$ is any positive integer. If ${\cal G}$ contains a maximal cycle
which generates a normal cyclic subgroup,
then ${\cal G}$ is a permutation matrix group.
\end{theorem}

Necessary preparations for the proof of the theorem are made in Section \ref{pre}.
For fundamentals of the group theory, please refer to \cite{AB,R}.
In Section \ref{cyclic case} we prove the theorem in the case
when the quotient group of ${\cal G}$ over the cyclic subgroup
generated by the maximal cycle is a cyclic group, see Theorem \ref{cyclic} below.
In Section \ref{non-cyclic case} we prove the theorem for the case
when the quotient group is not cyclic, see Theorem \ref{non-cyclic} below.

\section{Preparations}\label{pre}

We sketch necessary fundamentals from 
group theory, number theory and matrix theory, 
and formulate them in a way suitable to our later quotations.

The order of an element $C$ of a group ${\cal G}$ is denoted by
${\rm ord}_{\cal G}(C)$, or ${\rm ord}(C)$ in short
if the group is known from context. By $\langle C\rangle$
we denote the cyclic group generated by $C$.
Assume that ${\rm ord}(C)=d$. Then $\langle C\rangle$
is isomorphic to the additive group of the residue ring ${\Bbb Z}_{d}$
of the integer ring ${\Bbb Z}$ modulo $d$ by mapping
$j\in{\Bbb Z}_{d}$ to $C^j\in\langle C\rangle$.
In  this way, the properties of the cyclic group $\langle C\rangle$
are exactly corresponding to the properties of the additive group 
${\Bbb Z}_d$, some of which are listed as follows.
\begin{itemize}
\item
$\langle C^r\rangle=\langle C\rangle$ if and only if
$r\in{\Bbb Z}_d^*$, where ${\Bbb Z}_d^*$ denotes the multiplicative group
consisting of the reduced residue classes in ${\Bbb Z}_{d}$;
in that case $C^r$ is said to be a {\em generator}
of the cyclic group $\langle C\rangle$.
\item
Each automorphism $\alpha$ of $\langle C\rangle$
is corresponding to exactly one $r\in{\Bbb Z}_d^*$ such that
$\alpha(C^j)=C^{rj}$ for any $C^j\in\langle C\rangle$,
and $\mu_r$ is an automorphism of the additive group ${\Bbb Z}_d$,
where $\mu_r(j)=rj$ for any $j\in{\Bbb Z}_d$.
\item
If $\langle C\rangle$ is normal in ${\cal G}$,
then ${\cal G}/{\rm Cent}_{\cal G}(\langle C\rangle)$ is homomorphic to
a subgroup of ${\Bbb Z}_d^*$, where
${\rm Cent}_{\cal G}(\langle C\rangle)
 =\{A\in{\cal G}\mid A^{-1}CA=C\}$,
called the {\em centralizer} of $\langle C\rangle$ in $\cal G$.
If ${\rm Cent}_{\cal G}(\langle C\rangle)\subseteq\langle C\rangle$, then
$\langle C\rangle$ is said to be {\em self-centralized} in ${\cal G}$.
\end{itemize}

In this paper we are concerned with the cyclic group of order $2^n$
(with ${\Bbb Z}_{2^n}$ correspondingly).
The next lemma has appeared in \cite{DF15}.

\begin{lemma} \label{2-cyclic}
Assume that ${\cal G}=\langle C\rangle$ is a cyclic group of order
$|{\cal G}|=2^n$ 
and $0\le a<n$.
Let ${\cal G}^{2^a}=\{C^{2^a t}\mid 0\le t<2^{n-a}\}$. Then:
\begin{itemize}
\item[\rm(i)]  
${\cal G}^{2^a}=\langle C^{2^a}\rangle$
is a cyclic subgroup of ${\cal G}$ of order $2^{n-a}$ generated by $C^{2^a}$.
\item[\rm(ii)]
The mapping $C^j\mapsto C^{2^a j}$
is a homomorphism of ${\cal G}$ onto ${\cal G}^{2^a}$
with kernel ${\cal G}^{2^{n-a}}$.
Hence,  for any generator $G$ of the cyclic group ${\cal G}^{2^a}$
there are exactly $2^{a}$ generators of ${\cal G}$ which are mapped to $G$.
\qed
\end{itemize}
\end{lemma}

By $\Phi_{m}(x)$ we denote the cyclotomic polynomial of degree $m$, i.e.,
$\Phi_{m}(x)=\prod_{\theta}(x-\theta)$ with $\theta$ running over
the primitive $m$-th roots of unity.
Note that Lemma \ref{2-cyclic} (ii) (with $k=n-a$) implies that
\begin{equation}\label{cor 2-cyclic}
\prod_{\omega}(x^{2^a}-\omega)=\Phi_{2^{k+a}}(x)
\end{equation}
where the subscript $\omega$ of $\prod$
runs over the primitive $2^k$-th roots of unity.
In the specific case when $k=1$, since $-1$ is the unique primitive $2$-th root of unity,
we get a known formula: $\Phi_{2^{a+1}}(x)=x^{2^{a}}+1$.

\begin{remark}\label{Z^*}\rm
The structure of the multiplicative group ${\Bbb Z}_{2^n}^*$
is known (e.g., see \cite[Remark 9]{CDFL}). If $n\ge 3$ then the following hold:
\begin{itemize}
\item[(i)]
${\Bbb Z}_{2^n}^*=\langle 5\rangle\times\langle-1\rangle$
with~${\rm ord}_{{\Bbb Z}_{2^n}^*}(5)=2^{n-2}$ and
${\rm ord}_{{\Bbb Z}_{2^n}^*}(-1)=2$.
\item[(ii)]
$\langle 5\rangle=\{1\}\cup
\{1+2^{b}v\mid b=2,3,\cdots, n-1,~ v=1,3,\cdots, 2^{n-b}-1\}$,
and ${\rm ord}_{{\Bbb Z}_{2^n}^*}(1+2^{b}v)=2^{n-b}$.
\item[(iii)]
$-\langle 5\rangle=\{-1\}\cup
\{-1+2^{b}v\mid b=2,3,\cdots, n-1,~ v=1,3,\cdots, 2^{n-b}-1\}$,
and ${\rm ord}_{{\Bbb Z}_{2^n}^*}(-1+2^{b}v)=2^{n-b}$.
\end{itemize}
In our study, the behaviors of the elements $-1$,
$1+2^b v$ and $-1+2^b v$ are quite different.
But note that $(-1+2^b v)^\ell\in\langle 5\rangle$ provided $\ell$ is even.
\end{remark}

For the subgroups of ${\Bbb Z}_{2^n}^*$, we have the following result.

\begin{lemma}\label{sub Z^*}
If ${\cal H}$ is a subgroup of the multiplicative group
${\Bbb Z}_{2^n}^*$ where $n\ge 3$, then one of the following holds:
\begin{itemize}
\item[\rm(i)]
${\cal H}$ is a cyclic group, and, if the order $|{\cal H}|=2^a$,
then $0\le a\le n-2$ and ${\cal H}=\langle 1+2^{n-a}\rangle$,
or ${\cal H}=\langle -1+2^{n-a}\rangle$, or
${\cal H}=\langle -1\rangle$ (hence $a=1$).
\item[\rm(ii)]
${\cal H}=\langle -1\rangle\times\langle 1+2^{n-a}\rangle$
with $1\le a\le n-2$.
\end{itemize}
\end{lemma}

\pf By Remark~\ref{Z^*}, we have a surjective homomorphism
$\rho:{\Bbb Z}_{2^n}^*\to\langle 5\rangle$ such that $\rho(-1)=1$,
$\rho(1+2^{b}v)=1+2^b v$ and $\rho(-1+2^{b}v)=1+2^b(2^{n-b}-v)$.
Then $\rho({\cal H})$ is a subgroup of $\langle 5\rangle$,
hence $\rho({\cal H})$ is a cyclic group of order $2^a$ with $0\le a\le n-2$.
If $a=0$, then ${\cal H}=1$ or ${\cal H} =\langle -1\rangle$.
Assume that $a\ge 1$ in the following.
Take a $k\in{\cal H}$ such that $\rho(k)$ is a generator of $\rho({\cal H})$.
Then ${\cal H}=\langle k\rangle\times({\cal H}\cap\langle -1\rangle)$,
and $k=1+2^{n-a}$ or $k=-1+2^{n-a}$.
If ${\cal H}\cap\langle -1\rangle=1$, then (i) holds.
Otherwise, ${\cal H}\cap\langle -1\rangle=\langle -1\rangle$, hence
${\cal H}=\langle k\rangle\times\langle -1\rangle$;
if $k=-1+2^{n-a}$, we can replace $k$ by $-k$. In a word, (ii) holds.
\qed

\medskip
For any non-zero integer $k$,
by $\nu_2(k)$ we denote the $2$-adic valuation of $k$, i.e.,
$k=2^{\nu_2(k)}u$ for an odd integer $u$.

\begin{lemma}\label{1+r+} Let $r\in{\Bbb Z}_{2^n}^*$ such that
${\rm ord}_{{\Bbb Z}_{2^n}^*}(r)=2^a$ with $1\le a\le n-2$.
\begin{itemize}
\item[\rm(i)]
If $r=-1$, then $a=1$ and $r^{2^a-1}+r^{2^a-2}+\cdots+r+1 = 0$.
\item[\rm(ii)]
If $r=1+2^{n-a}v$ with $v$ being odd, then
$$
\nu_2(r^{2^a-1}+r^{2^a-2}+\cdots+r+1)=a.
$$
\item[(iii)]
If $r=-1+2^{n-a}v$ with $v$ being odd, then
$$
\nu_2(r^{2^a-1}+r^{2^a-2}+\cdots+r+1)=n-1.
$$
\end{itemize}
\end{lemma}

\pf (i) is obvious.

(ii).~ Since $\nu_2(2^a-j)=\nu_2(j)$ for $0< j< 2^a$,
for the binomial coefficients $\binom{2^a}{k}$ we get that
\begin{equation}\label{binom coef}
\nu_2\binom{2^a}{k}=\nu_2
\Big(2^a\cdot\frac{2^a-1}{1}\cdot\cdots\frac{2^a-(k-1)}{k-1}\cdot\frac{1}{k}\Big)
=a-\nu_2(k).
\end{equation}
So
$$
r^{2^a-1}+\cdots+r+1=\frac{r^{2^a}-1}{r-1}
=\frac{(1+2^{n-a}v)^{2^a}-1}{2^{n-a}v}
\equiv 2^a\pmod{2^{a+1}}.
$$

(iii).~ Noting that $r-1=2u$ with $u$ being odd,
by Eqn \eqref{binom coef} we obtain that
$$
r^{2^a-1}+\cdots+r+1=\frac{r^{2^a}-1}{r-1}
=\frac{(-1+2^{n-a}v)^{2^a}-1}{2u}\equiv 2^{n-1}v'\pmod{2^n};
$$
where $v'\equiv -vu^{-1}\!\pmod{2^n}$.
\qed

\begin{lemma} \label{split}
Let ${\cal G}=\langle A,C\rangle$
be a finite group generated by $A$ and $C$.
If $\langle C\rangle$ is a self-centralized 
normal cyclic subgroup of order $2^n$
and $A^{-1}CA=C^r$ for an $r\in{\Bbb Z}_{2^n}$
with ${\rm ord}_{{\Bbb Z}_{2^n}^*}(r)=2^a$,
then $|{\cal G}|=2^{n+a}$ and:
\begin{itemize}
\item[\rm(i)] 
If  $r\not\equiv -1~({\rm mod}~2^n)$,  then
there is an $A'\in{\cal G}$ such that ${\cal G}=\langle A',C\rangle$ with
$A'^{-1}CA'=C^r$ and $A'^{2^a}=1$.
\item[\rm(ii)]
If $r\equiv -1~({\rm mod}~2^n)$ (hence $a=1$),
then one of the following two holds:
\begin{itemize}
\item[\rm(ii.1)] 
$A^2=1$; in that case $(AC^k)^2=1$ for any $C^k\in\langle C\rangle$
(i.e., $\cal G$ is the dihedral group of order $2^{n+1}$).
\item[\rm(ii.2)]
$A^2=C^{2^{n-1}}$; in that case $(AC^k)^2=C^{2^{n-1}}$
for any $C^k\in\langle C\rangle$
(i.e., $\cal G$ is the generalized quaternion group of order $2^{n+1}$).
\end{itemize}
\end{itemize}
\end{lemma}

\pf
Since $\langle C\rangle$ is self-centralized,
the quotient group ${\cal G}/\langle C\rangle
\cong\langle A\rangle/\langle A\rangle\cap\langle C\rangle$
is isomorphic to a cyclic subgroup $\langle r\rangle$ of the group ${\Bbb Z}_{2^n}^*$.
So $|\langle A\rangle/\langle A\rangle\cap\langle C\rangle|=2^a$,
 $|{\cal G}|=2^{n+a}$ and $A^{2^a}\in\langle C\rangle$.
Let $\langle C\rangle^A=\{C^t\mid A^{-1}C^tA=C^t\}$,
i.e., the centralizer in $\langle C\rangle$ of $A$.
Then $A^{2^a}\in\langle C\rangle^A$.
Since $A^{-1}C^tA=C^{tr}$,
\begin{equation}\label{C^A}
 A^{-1}C^tA=C^t
 ~~~{\rm iff}~~~ C^{tr-t}=1~~~{\rm iff}~~~
  t(r-1)\equiv 0\hskip-5pt\pmod{2^n}.
\end{equation}
Suppose that $r=-1+2^{n-a}v$ where $v$ is zero (i.e., $r=-1$) or odd.
Since $n-a\ge 2$,  we have $r-1=-2+v2^{n-a}=2u$ with $u$ being odd.
$$
 A^{-1}C^tA=C^t ~~~{\rm iff}~~~ 2ut\equiv 0\hskip-5pt\pmod{2^n}
 ~~~{\rm iff}~~~ t\equiv 0\hskip-5pt\pmod{2^{n-1}}.
$$
Thus $\langle C\rangle^A =\langle C^{2^{n-1}}\rangle$
which is a group of order $2$ consisting of two elements $\{1,C^{2^{n-1}}\}$.
So
$$A^{2^a}\in\langle C^{2^{n-1}}\rangle=\{1,C^{2^{n-1}}\}.$$
In particular, if $v=0$, i.e., $r=-1$, then we obtain the conclusion (ii).

Next, we assume that $r=-1+2^{n-a}v$ with $v$ being odd.
Then $A^{2^a}=1$ or $A^{2^a}=C^{2^{n-1}}$.
If $A^{2^a}=1$ then we are done. Suppose that $A^{2^a}=C^{2^{n-1}}$.
Note that for any integer $k$ and any positive integer $j$, by a direct computation
we have the following formula:
\begin{equation}\label{(AC)^j}
(AC^k)^j=(AC^k)\cdots(AC^k)(AC^k)=A^{j}C^{k(r^{j-1}+\cdots+r+1)}.
\end{equation}
By Lemma \ref{1+r+}(iii),
$r^{2^a-1}+\cdots+r+1=2^{n-1}u$ with $u$ being odd.
We have an integer $u'$ such that $u'u\equiv 1~({\rm mod}~2^n)$.
Taking $j=2^a$ and $k=u'$ in Eqn \eqref{(AC)^j},  we obtain
$$
(AC^{u'})^{2^a}=A^{2^a}C^{u'(r^{2^a-1}+\cdots+r+1)}
=A^{2^a}C^{2^{n-1}u'u}=C^{2^{n-1}}C^{2^{n-1}}=1.
$$
Replacing $A$ by $A'=AC^{u'}$, we get that
${\cal G}=\langle A',C\rangle$, $A'^{-1}CA'=C^{r}$ and $A'^{2^a}=1$.
So (i) holds.

Finally, assume that $r=1+2^{n-a}v$ with $v$ being odd.  Then
$$
t(r-1)\equiv 0\hskip-5pt\pmod{2^n}~~~{\rm iff}~~~
2^{n-a}vt \equiv 0\hskip-5pt\pmod{2^n}~~~{\rm iff}~~~
 t \equiv 0\hskip-5pt\pmod{2^{a}}.
$$
By Eqn \eqref{C^A}, $\langle C\rangle ^A=\langle C^{2^a}\rangle$;
hence
$A^{2^a}\in\langle C\rangle^A=\langle C^{2^a}\rangle.$
So, we can find an integer $h$ such that $A^{2^a}C^{2^ah}=1$.
By Lemma \ref{1+r+}(ii) we can assume that
$r^{2^a-1}+\cdots+r+1=2^a u$ with $u$ being odd.
We have an integer $u'$ such that $u'u\equiv 1~({\rm mod}~2^n)$,
Taking $j=2^a$ and $k=u'h$ in Eqn \eqref{(AC)^j},  we obtain
$$
(AC^{u'h})^{2^a}=A^{2^a}C^{u'h(r^{2^a-1}+\cdots+r+1)}=A^{2^a}C^{2^auu'h}=1.
$$
Replacing $A$ by $A'=AC^{u'h}$, we get the conclusion (i).
\qed

\medskip
The next lemma shows that the condition ``self-centralized'' is 
usually satisfied in our study.

\begin{lemma}[{\cite[Proposition 4.2]{C07}}]\label{self-centralized}
If ${\cal G}$ is a permutation-like matrix group and
$C\in{\cal G}$ is a maximal cycle, then $\langle C\rangle$
is self-centralized in ${\cal G}$.
\end{lemma}

There is also a proof of this lemma in \cite[Lemma 2.7]{DF15}.

\medskip
For any complex matrix $A$, by ${\rm char}_A(x)$ we denote
the characteristic polynomial of $A$, i.e.,
${\rm char}_A(x)=\det(xI-A)$, where $I$ denotes the identity matrix.

\begin{lemma}[{\cite[Lemma 2.1]{DF}, \cite[Lemma 2.3]{DF15}}]
\label{permutable}
The following two are equivalent to each other:
\begin{itemize}
\item[\rm(i)]
$A$ is similar to a permutation matrix.
\item[\rm(ii)]
$A$ is diagonalizable and ${\rm char}_A(x)=\prod_{i}(x^{\ell_i}-1)$.
\end{itemize}
If it is the case, then each factor $x^{\ell_i}-1$ of ${\rm char}_A(x)$
corresponds to exactly one $\ell_i$-cycle of the cycle decomposition
of the permutation of the permutation matrix. \qed
\end{lemma}

Note that $x^n-1=\prod_{k|n}\Phi_k(x)$, where $\Phi_k(x)$ is the
cyclotomic polynomial of degree $k$, and ``$k|n$'' stands for
that $k$ divides $n$.  
We get an immediate consequence of the above lemma as follows.

\begin{corollary}\label{cor permutable}
Let $A$ be a matrix similar to a permutation matrix, and
$m,n$ be positive integers. If ${\Phi_n(x)}^m\,\big|\,{\rm char}_A(x)$,
then $\Phi_k(x)^m\,\big|\,{\rm char}_A(x)$ for any $k\,|\,n$.
\qed
\end{corollary}

\section{The quotient ${\cal G}/\langle C\rangle$ is cyclic} \label{cyclic case}

In this section we always assume:
\begin{itemize}
\item 
$A, C\in{\rm GL}_{2^n}({\Bbb C})$;
\item 
$C$ is a maximal cycle;
\item 
$A^{-1}CA=C^r$, where $r\in{\Bbb Z}_{2^n}^*$
and ${\rm ord}_{{\Bbb Z}_{2^n}^*}(r)=2^a$ with $0\le a\le n-2$.
\end{itemize}
Let $\lambda$ be a primitive $2^n$-th root of unity. Then
\begin{equation}\label{spectrum}
\{\lambda^j\mid j\in{\Bbb Z}_{2^n}\}
  =\{1=\lambda^0,\lambda, \lambda^2, \cdots,\lambda^{2^n-1}\}
\end{equation}
is the spectrum (i.e., the set of eigenvalues) of $C$.
The eigen-subspace of each eigenvalue $\lambda^j$ of $C$,
denoted by ${\rm E}(\lambda^j)$, has dimension $1$; and the complex space
${\Bbb C}^{2^n}$ is decomposed into a direct sum of the eigen-subspaces:
\begin{equation}\label{eigen-subspace}
{\Bbb C}^{2^n}=\bigoplus_{j\in{\Bbb Z}_{2^n}}{\rm E}(\lambda^j)=
{\rm E}(\lambda^0)\oplus{\rm E}(\lambda)
\oplus {\rm E}(\lambda^2)\oplus\cdots\oplus {\rm E}(\lambda^{2^n-1}).
\end{equation}
Taking any non-zero vector $e_j\in{\rm E}(\lambda^j)$ 
for every index $j\in{\Bbb Z}_{2^n}$,
we get a basis of ${\Bbb C}^{2^n}$:
\begin{equation}\label{basis}
{\cal E}=\{e_j\mid j\in{\Bbb Z}_{2^n}\}=\{e_0,e_1,e_2,\cdots,e_{2^n-1}\},
\end{equation}
which will play a key role in our study.
For convenience we denote the diagonal blocked matrix
$\begin{pmatrix}B_1\\ & \ddots\\ &&B_m\end{pmatrix}$
by $B_1\oplus\cdots\oplus B_m$.
With the basis ${\cal E}$ in Eqn~\eqref{basis},
the transformation $C$ on ${\Bbb C}^{2^n}$ has the diagonal  matrix
$$
 C|_{\cal E}=1\oplus \lambda\oplus \lambda^2\oplus\cdots\oplus \lambda^{2^n-1}.
$$
For the above, please see \cite[Lemma 2.4]{DF15}.

Consider the matrix $A|_{\cal E}$ with respect to 
the basis ${\cal E}$ in Eqn \eqref{basis}.
Recall that $A^{-1}CA=C^r$ where $r\in{\Bbb Z}_{2^n}^*$ with
${\rm ord}_{{\Bbb Z}_{2^n}^*}(r)=2^a$.  It is easy to see that
\begin{equation}\label{permute E}
A{\rm E}(\lambda^{j})={\rm E}(\lambda^{rj}),\qquad
\forall ~ j\in{\Bbb Z}_{2^n};
\end{equation}
i.e., $A$ permutes the eigen-subspaces
${\rm E}(\lambda^j)$ for $j\in{\Bbb Z}_{2^n}$ the same as
$\mu_r$ permutes ${\Bbb Z}_{2^n}$, cf. \cite[Lemma~2.5]{DF15}.
For any $\mu_r$-orbit $\Gamma_k=\{j_k,rj_k,\cdots,r^{2^b-1}j_k\}$
on ${\Bbb Z}_{2^n}$ of length $2^b$ ($b\le a$ of course), the subset
${\cal E}_k=\{e_{j_k}, e_{rj_k},\cdots,e_{r^{2^b-1}j_k}\}$ of ${\cal E}$
is a basis of the subspace
$V_k=\bigoplus_{i=0}^{2^b-1}{\rm E}(\lambda^{r^ij_k})$.
And there are $\alpha_1,\cdots,\alpha_{2^b}\in{\Bbb C}$ such that
$Ae_{j_k}=\alpha_1 e_{rj_k}$, $Ae_{rj_k}=\alpha_2 e_{r^2j_k}$, $\cdots$,
$Ae_{r^{2^b-1}j_k}=\alpha_{2^b} e_{j_k}$.
That is, $V_k$ is invariant by $A$,
and the matrix of $A$ with respect to the basis ${\cal E}_k$ is
\begin{equation}\label{monomial cycle}
   A|_{{\cal E}_k}= \begin{pmatrix}0&\cdots&0&\alpha_{2^b}\\
    \alpha_1 &\ddots&\ddots&0\\ \vdots&\ddots&\ddots&\vdots\\
     0&\cdots&\alpha_{2^b-1}&0\end{pmatrix}_{2^b\times 2^b}.
\end{equation}
If ${\rm ord}(A)=2^{a'}$ ($a'\ge a$ of course), then
$\omega_{2^b}=\alpha_1\cdots\alpha_{2^b}$ being a $2^{a'-b}$-th root of unity.
Replacing ${\cal E}_k$ by
\begin{equation}\label{cycle basis}
{\cal E}_k=\{e_{j_k},Ae_{j_k},\cdots,A^{2^b-1}e_{j_k}\},
\end{equation}
we get
\begin{equation}\label{1-monomial cycle}
   A|_{{\cal E}_k}= \begin{pmatrix}0&\cdots&0&\omega_{2^b}\\ 1&\ddots&\ddots&0\\
    \vdots&\ddots&\ddots&\vdots\\ 0&\cdots&1&0\end{pmatrix}_{2^b\times 2^b}.
\end{equation}
In particular,
$A|_{{\cal E}_k}$ is a permutation matrix once $a'=b$
(because $\omega_{2^b}=1$ once $a'=b$);
in that case, $A$ permutes the set ${\cal E}_k$ cyclically.

Similarly to $\langle C\rangle$, by $\langle A,C\rangle$ we denote
the group generated by $A$ and $C$.

\begin{lemma}[{\cite[Proposition 2.3]{DF},
 \cite[Proposition 2.6]{DF15}}]\label{SC}
Let notations be as above.
If there is a choice of the basis ${\cal E}$ in Eqn \eqref{basis} such that
the matrix $A|_{\cal E}$ is a permutation matrix,
then the matrix group $\langle A,C\rangle$ 
is a permutation matrix group. 
\end{lemma}

The statement ``the matrix $A|_{\cal E}$ is a permutation matrix''
is equivalent to saying that ``the set ${\cal E}$ is $A$-stable'',
i.e., $A$ permutes the elements of the set~${\cal E}$.
In that case we call ${\cal E}$ an {\em $A$-permutation basis}.
The rest of this section contributes mainly to   
finding such a basis ${\cal E}$.

We partition ${\Bbb Z}_{2^n}$ into a disjoint union
${\Bbb Z}_{2^n}=2{\Bbb Z}_{2^n}\bigcup {\Bbb Z}_{2^n}^*$, where
$$\begin{array}{l}
2{\Bbb Z}_{2^n}=\{2k~({\rm mod}~2^n)\mid k\in{\Bbb Z}_{2^n}\}
 =\{0,2,\cdots,2(2^{n-1}-1)\},\\[3mm]
{\Bbb Z}_{2^n}^*=1+2{\Bbb Z}_{2^n}
=\{1+2k~({\rm mod}~2^n)\mid k\in{\Bbb Z}_{2^n}\}
=\{1,3,\cdots,2^{n}-1\}.
\end{array}$$
Correspondingly, the basis ${\cal E}$ in Eqn \eqref{basis} is partitioned into
${\cal E}={\cal E}^2\cup{\cal E}^*$ where
${\cal E}^2=\{e_j\,|\, j\in 2{\Bbb Z}_{2^n}\}$
and ${\cal E}^*=\{e_j\,|\, j\in {\Bbb Z}_{2^n}^*\}$.
And by Eqn \eqref{eigen-subspace}, ${\Bbb C}^{2^n}$ is decomposed into
a direct sum of two subspaces: ${\Bbb C}^{2^n}=V^2\oplus V^*$ where
\begin{equation}\label{V2V*}
\begin{array}{l}
V^{2} =\bigoplus\limits_{j\in 2{\Bbb Z}_{2^n}}{\rm E}(\lambda^j)
=\bigoplus\limits_{e\in{\cal E}^2}{\Bbb C}e,\\[5mm]
 V^{*} = \bigoplus\limits_{j\in {\Bbb Z}_{2^n}^*}{\rm E}(\lambda^j)
=\bigoplus\limits_{e\in{\cal E}^*}{\Bbb C}e.
\end{array}
\end{equation}
Both $2{\Bbb Z}_{2^n}$ and ${\Bbb Z}_{2^n}^*$ are $\mu_r$-invariant.
So, for any $k$ with $0\le k<2^n$, both $V^{2}$ and $V^{*}$
are $AC^k$-invariant subspaces of ${\Bbb C}^{2^n}$.
By $AC^k|_{V^2}$ and $AC^k|_{V^*}$ we denote the linear transformations
of $AC^k$ restricted to $V^2$ and $V^*$ respectively.
Correspondingly,
$AC^k|_{{\cal E}^2}$ and $AC^k|_{{\cal E}^*}$ are
matrices of $AC^k|_{V^2}$ and $AC^k|_{V^*}$ respectively.

\begin{remark}\label{e_0}\rm
The two elements $0,2^{n-1}\in{\Bbb Z}_{2^n}$ are always fixed by
$\mu_r$ for any $r\in{\Bbb Z}_{2^n}^*$.
So, the two eigen-subspaces ${\rm E}(\lambda^{0})={\rm E}(1)$ and
${\rm E}(\lambda^{2^{n-1}})={\rm E}(-1)$
are always $A$-invariant subspaces.

(i).~ Note that:
{\it if $\langle A,C\rangle$ is a
permutation-like matrix group, then $Ae_{0}=e_{0}$
for any $e_0\in{\rm E}(\lambda^0)$.}~
Because: in that case $\sum_{G\in\langle A,C\rangle}G\ne 0$
(since the traces of permutation matrices are non-negative integers,
and the trace of $I$ is $2^n>0$),  there is a
vector $f\in{\Bbb C}^{2^n}$ such that 
$f_0=\sum_{G\in\langle A,C\rangle}Gf\ne 0$;
it is easy to see that $Cf_0=f_0$ and $Af_0=f_0$;
hence $f_0\in {\rm E}(\lambda^0)$, and
$Ae_{0}=e_{0}$ for any $e_0\in{\rm E}(\lambda^0)$.

(ii). It is obvious that, {\it if $A|_{\cal E}$ is a permutation matrix,
then $Ae_{2^{n-1}}=e_{2^{n-1}}$.}
\end{remark}

\begin{lemma}\label{V^*}
Assume that $A, C\in{\rm GL}_{2^n}({\Bbb C})$, 
$C$ is a maximal cycle and
$A^{-1}CA=C^r$ where $r\in{\Bbb Z}_{2^n}^*$
with ${\rm ord}_{{\Bbb Z}_{2^n}^*}(r)=2^a$.
If $A^{2^a}=I$,
then there is a choice of the basis ${\cal E}^*$ of
$V^*$ in Eqn \eqref{V2V*} such that
$A|_{{\cal E}^*}$ is a permutation matrix.
\end{lemma}

\pf Since ${\rm ord}_{{\Bbb Z}_{2^n}^*}(r)=2^a$,
the length of any $\mu_r$-orbit on ${\Bbb Z}_{2^n}^*$ equals to $2^a$.
Note that ${\rm ord}(A)=2^a$.
The conclusion follows from Eqn \eqref{1-monomial cycle}.
\qed

\medskip
The following lemma exhibits the information for the case when $r=-1$,
which is enough for our later quotations.

\begin{lemma} \label{BC^k}
Assume that $A, C\in{\rm GL}_{2^n}({\Bbb C})$,
 $C$ is a maximal cycle and $A^{-1}CA=C^{-1}$.
If ${\cal G}=\langle A,C\rangle$ is a
permutation-like matrix group, then $A^2=I$
(i.e., ${\cal G}$ is a dihedral group of order $2^{n+1}$) and
one of the following two holds:
\begin{itemize}
\item[\rm(i)]
$Ae_{2^{n-1}}=e_{2^{n-1}}$ for
$e_{2^{n-1}}\in{\rm E}(\lambda^{2^{n-1}})$, and
$${\rm char}_{AC^{k}}(x)=
 \begin{cases} (x-1)^2(x^2-1)^{2^{n-1}-1}, &  \mbox{$k$ is even};\\
  (x^2-1)^{2^{n-1}}, &  \mbox{$k$ is odd}. \end{cases}
$$
In that case, there is a choice of the basis ${\cal E}$ in Eqn \eqref{basis}
such that $A|_{\cal E}$ is a permutation matrix.
\item[\rm(ii)]
$Ae_{2^{n-1}}=-e_{2^{n-1}}$ for
$e_{2^{n-1}}\in{\rm E}(\lambda^{2^{n-1}})$, and
$${\rm char}_{AC^{k}}(x)=
 \begin{cases} (x^2-1)^{2^{n-1}}, &  \mbox{$k$ is even};\\
   (x-1)^2(x^2-1)^{2^{n-1}-1}, &  \mbox{$k$ is odd}. \end{cases}
$$
\end{itemize}
\end{lemma}

\pf It is easy to see that ${\Bbb Z}_{2^n}$ is partitioned
in to $\mu_{-1}$-orbits as follows:
\begin{equation*}
\begin{array}{c}
2{\Bbb Z}_{2^n}=
\{0\}\cup\{2^{n-1}\}\cup\{2,-2\}\cup\cdots\cup\{2^{n-1}-2,\;-2^{n-1}+2\},
\\[7pt]
{\Bbb Z}_{2^n}^*=
\{1,-1\}\cup\{3,-3\}\cup\cdots\cup\{2^{n-1}-1,\;-2^{n-1}+1\}.
\end{array}
\end{equation*}
Correspondingly, $C$ is diagonally blocked as follows:
\begin{align*}
C|_{\cal E}&=C|_{{\cal E}^2}\oplus C|_{{\cal E}^*}\\
&=(1\oplus -1\oplus \lambda^2\oplus \lambda^{-2}\oplus\cdots)\oplus
 (\lambda^1\oplus \lambda^{-1}\oplus \lambda^3\oplus \lambda^{-3}\oplus\cdots).
\end{align*}

By Lemma~\ref{split}, either $A^2=I$ or $A^2=C^{2^{n-1}}$.
We prove that $A^2\ne C^{2^{n-1}}$ by contradiction.
Suppose that $A^2=C^{2^{n-1}}$, i.e.,
$A|_{\cal E}=A|_{{\cal E}^*}\oplus A|_{{\cal E}^2}$ with
$$
A^2|_{{\cal E}^2}=I_{2^{n-1}}, \qquad A^2|_{{\cal E}^*}=-I_{2^{n-1}}.
$$
By Eqn \eqref{1-monomial cycle},
After a suitable choice of the basis ${\cal E}$ we can assume that
$$
A|_{{\cal E}^2}=1\oplus \pm1\oplus
\overbrace{P\oplus P\oplus\cdots\oplus P}^{2^{n-2}-1}
\qquad \mbox{where}~~
P=\begin{pmatrix} 0& 1\\1&0 \end{pmatrix},
$$
and
$$
A|_{{\cal E}^*}=\overbrace{Q\oplus Q\oplus \cdots\oplus Q}^{2^{n-2}}
\qquad \mbox{where}~~
Q=\begin{pmatrix} 0& -1\\1&0 \end{pmatrix}.
$$
So, ${\rm char}_{A|_{{\cal E}^2}}(x)=(x-1)(x\mp 1)(x^2-1)^{2^{n-2}-1}$
and ${\rm char}_{A|_{{\cal E}^*}}(x)=(x^2+1)^{2^{n-2}}$.
We get
$$
{\rm char}_A(x)=(x-1)(x\mp 1)(x^2+1)(x^4-1)^{2^{n-2}-1}.
$$
By Lemma~\ref{permutable}, the second factor of the right hand side
has to be $x+1$. Hence
$$
 A|_{{\cal E}}=A|_{{\cal E}^2}\oplus A|_{{\cal E}^*}=
(1\oplus -1\oplus P\oplus P\oplus\cdots)\oplus(Q\oplus Q\oplus \cdots).
$$
Then we get the matrix of $AC$ as follows:
$$
AC=\Bigg(1\oplus 1\oplus
\begin{pmatrix} & \lambda^2\\\lambda^{-2} \end{pmatrix}\oplus\cdots\Bigg)
\oplus \Bigg(\begin{pmatrix} & -\lambda^{-1}\\ \lambda \end{pmatrix}\oplus
\begin{pmatrix} & -\lambda^{-3}\\ \lambda^{3} \end{pmatrix}\oplus\cdots\Bigg)
$$
and the characteristic polynomial as follows:
\begin{eqnarray*}
{\rm char}_{AC}(x)&=&(x-1)^2(x^2-1)^{2^{n-2}-1}(x^2+1)^{2^{n-2}}\\
&=&(x-1)^2(x^2+1)(x^4-1)^{2^{n-2}-1}.
\end{eqnarray*}
By Lemma \ref{permutable}, $AC$ is not similar to a permutation matrix.
However, the assumption of the lemma says that
$AC$ is similar to a permutation matrix. The contradiction
forces that $A^2\ne C^{2^{n-1}}$.

Thus, we obtain that $A^2=I$ and
$$A|_{{\cal E}^*}=\overbrace{P\oplus P\oplus \cdots\oplus P}^{2^{n-2}},
 \qquad
A|_{{\cal E}}=1\oplus \pm1\oplus
\overbrace{P\oplus P\oplus\cdots\oplus P}^{2^{n-2}-1}.$$
It is easy to check that:

$\bullet$~ if
$A|_{{\cal E}}=1\oplus 1\oplus
\overbrace{P\oplus\cdots\oplus P}^{2^{n-1}-1}$,
then (i) holds;

$\bullet$~ otherwise, $A|_{{\cal E}}=1\oplus -1\oplus
\overbrace{P\oplus\cdots\oplus P}^{2^{n-1}-1}$
and (ii) holds.
\qed

\medskip
In the following we consider the case when $r\not\equiv-1\!\pmod{2^n}$ and
find a basis ${\cal E}$ in Eqn \eqref{basis}
such that $A|_{\cal E}$ is a permutation matrix.
We start it with an investigation of
the characteristic polynomial of $AC^k$.

\begin{lemma} \label{AC^k}
Let $A, C\in{\rm GL}_{2^n}({\Bbb C})$.
Assume that $C$ is a maximal cycle,
$A^{-1}CA=C^r$ where $-1\ne r\in{\Bbb Z}_{2^n}^*$
with ${\rm ord}_{{\Bbb Z}_{2^{n}}^*}(r)=2^{a}$,
and $A^{2^{a}}=I$.
Let $V^2$ and $V^*$ be the subspaces in Eqn \eqref{V2V*}.
\begin{itemize}
\item[\rm(i)]
If $r=1+2^{n-a} v$ where $0\le a\le n-2$ and $v$ is odd, then
$${\rm char}_{(AC^k)|_{V^*}}(x)
=\begin{cases}
 {\Phi_{2^{n-\nu_2(k)}}(x)}^{2^{\nu_2(k)}}, & 0\le\nu_2(k)<n-a;\\[7pt]
 {(x^{2^a}-1)}^{2^{n-a-1}}, &  \nu_2(k)\ge n-a,~{\rm or}~ k=0. \end{cases}
$$
\item[\rm(ii)]
If $r=-1+2^{n-a} v$ where $1\le a\le n-2$ and $v$ is odd, then
$${\rm char}_{(AC^k)|_{V^*}}(x)
=\begin{cases}{(x^{2^a}+1)}^{2^{n-a-1}}=
{\Phi_{2^{a+1}}(x)}^{2^{n-a-1}}, & \mbox{$k$ is odd};\\[7pt]
 {(x^{2^a}-1)}^{2^{n-a-1}}, &  \mbox{$k$ is even}. \end{cases}
$$
\end{itemize}
\end{lemma}

\pf
By the assumption, ${\Bbb Z}_{2^n}^*$  is partitioned into $\mu_r$-orbits
$\Gamma_1$, $\cdots$, $\Gamma_h$ with the lengths all being $2^{a}$
and the number $h=2^{n-a-1}$:
$$
\Gamma_1=\{j_1,rj_1,\cdots,r^{2^{a}-1}j_1\},~  \cdots,~
\Gamma_{h}=\{j_h, rj_h, \cdots, r^{2^{a}-1}j_h\}.
$$
By Eqn \eqref{basis} and Eqn \eqref{cycle basis},
we can assume that
 ${\cal E}^*={\cal E}^*_1\cup\cdots\cup{\cal E}^*_h$ with
\begin{equation}\label{E_i}
{\cal E}^*_i=\big\{e_{j_i},~Ae_{j_i},~\cdots,~A^{r^{2^a-1}}e_{j_i}\big\},
\qquad i=1,\cdots,h.
\end{equation}
Denote the matrices of $A$ and $C$ restricted to ${\cal E}^*_i$ by
$A_i=A|_{{\cal E}^*_i}$ and $C_i=C|_{{\cal E}^*_i}$ for $i=1,\cdots,h$.
Since $A_i^{2^a}=I$, by Lemma \ref{V^*} we have
\begin{equation*}
A_{i}=
   \begin{pmatrix}0&\cdots&0&1\\ 1&\ddots&\ddots&0\\
          \vdots&\ddots&\ddots&\vdots\\
          0&\cdots&1&0   \end{pmatrix}_{2^{a}\times 2^{a}},~
 C_{i}^k=
   \begin{pmatrix}\lambda^{kj_i} \\ &\lambda^{kj_i r}\\
   &&\ddots\\ &&& \lambda^{kj_i r^{2^a-1}}
 \end{pmatrix}_{2^{a}\times 2^{a}}.
\end{equation*}
Thus
$
AC^k|_{{\cal E}^*}=\bigoplus_{i=1}^h A_iC_i^k,
$
and
$
{\rm char}_{A_iC_i^k}(x)=x^{2^a}-\lambda^{j_ik(1+r+\cdots+r^{2^a-1})}.
$
The conclusion is obviously true if $k=0$. So we further assume that $k\ne 0$.

Writing $k=k'2^{\nu_2(k)}$ with $k'$ being odd, and
writing $1+r+\cdots+r^{2^a-1}=a'2^{a''}$ with $a'$ being odd and
$a''=\nu_2(1+r+\cdots+r^{2^a-1})$, we have
\begin{equation*}
{\rm char}_{A_iC_i^k}(x)=x^{2^a}-\lambda^{j_ia'k'2^{a''+\nu_2(k)}},\qquad
i=1,\cdots,h.
\end{equation*}
If $a''+\nu_2(k)\ge n$, then $\lambda^{j_ia'k'2^{a''+\nu_2(k)}}=1$.
However, if $a''+\nu_2(k)<n$, then $\lambda^{j_ia'k'2^{a''+\nu_2(k)}}$ is a
primitive $2^{n-a''-\nu_2(k)}$-th root of unity,
and, by Lemma~\ref{2-cyclic}~(ii), the collection of
$$\lambda^{j_i a'k' 2^{a''+\nu_2(k)}},\qquad  i=1,\cdots, h,$$
is just the collection of all primitive $2^{n-a''-\nu_2(k)}$-th roots of unity,
each of which appears with multiplicity
$\frac{h}{2^{{n-a''-\nu_2(k)}-1}}=2^{a''-a+\nu_2(k)}.$
By Eqn \eqref{cor 2-cyclic}, we obtain that
\begin{equation*}
{\rm char}_{(AC^k)|_{V^*}}(x)=
\begin{cases} (x^{2^a}-1)^{n-a-1}, & a''+\nu_2(k)\ge n;\\
{\Phi_{2^{n+a-a''-\nu_2(k)}}(x)}^{2^{a''-a+\nu_2(k)}}, & a''+\nu_2(k)<n.
\end{cases}
\end{equation*}

If $r=1+2^{n-a} v$ where $0\le a\le n-2$ and $v$ is odd,
then, by Lemma \ref{1+r+}(ii), $a''=a$; hence (i) is proved.

Next, assume that $r=-1+2^{n-a} v$ where $1\le a\le n-2$ and $v$ is odd.
By Lemma \ref{1+r+}(iii), $a''=n-1$. So $a''+\nu_2(k)< n$ if and only if
$k$ is odd; in that case, $n+a-a''-\nu_2(k)=a+1$.
Noting that $\Phi_{2^{a+1}}(x)=x^{2^a}+1$,
we are done for the conclusion (ii).
\qed

\medskip
The second step is a preparation for induction.
Let $r\in{\Bbb Z}_{2^n}^*$. Then $\mu_r$ (the action by multiplication) acts on both
${\Bbb Z}_{2^{n-1}}$ and $2{\Bbb Z}_{2^{n}}$, and the map
\begin{equation}\label{eq 2Z}
 {\Bbb Z}_{2^{n-1}} \longrightarrow 2{\Bbb Z}_{2^{n}},~~
  t\longmapsto 2t,
\end{equation}
is an isomorphism between the $\mu_r$-acted sets
${\Bbb Z}_{2^{n-1}}$ and $2{\Bbb Z}_{2^{n}}$.

\begin{lemma}\label{V^2}
Let notation be as in Lemma \ref{AC^k}. Assume further that
the matrix group $\langle A,C\rangle$ is permutation-like.
Then $\big\langle A|_{V^2},C|_{V^2}\big\rangle$
is a permutation-like matrix group of dimension $2^{n-1}$
with $C|_{V^2}$ being a maximal cycle.
\end{lemma}

\pf Since $C|_{{\cal E}^2}=\bigoplus_{j=0}^{2^{n-1}-1}\lambda^j$,
with the isomorphism \eqref{eq 2Z},
$C|_{V^2}$ is a maximal cycle of length $2^{n-1}$.
Note that
$\langle A,C\rangle=\{A^\ell C^k\mid 0\le\ell<2^a,\;0\le k<2^n\}.$
To complete the proof, it is enough to
show that any $(A^\ell C^k)|_{V^2}$ is similar to a permutation matrix.

Let $A'=A^\ell$, $a'=a-\nu_2(\ell)$ and $r'=r^\ell$.
Then $A'^{-1}CA'=C^{r'}$ and
${\rm ord}(A')={\rm ord}_{{\Bbb Z}_{2^n}^*}(r')=2^{a'}$.
We compute
${\rm char}_{(A^\ell C^k)|_{V^2}}(x)={\rm char}_{(A' C^k)|_{V^2}}(x)$
in two cases.

{\em Case 1}: $r\in\langle 5\rangle$ or $\nu_2(\ell)\ge 1$
(cf. Remark \ref{Z^*}). We can use Lemma~\ref{AC^k}(i) to
$A'C^k$ and get that
\begin{equation*}
{\rm char}_{(A^\ell C^k)|_{V^*}}(x)
=\begin{cases}{\Phi_{2^{n-\nu_2(k)}}(x)}^{2^{\nu_2(k)}},
  & 0\le \nu_2(k)<n-a';\\[7pt]
 {(x^{2^{a'}}-1)}^{2^{n-a'-1}}, &  \nu_2(k)\ge n-a',~{\rm or}~k=0 . \end{cases}
\end{equation*}
Since $A^\ell C^k$ is similar to a permutation matrix,
by Lemma \ref{permutable} and its corollary,
\begin{equation}\label{r=5 V^2}
{\rm char}_{(A^\ell C^k)|_{V^2}}(x)=\begin{cases}
{(x^{2^{n-\nu_2(k)-1}}-1)}^{2^{\nu_2(k)}}, & 0\le \nu_2(k)<n-a';\\[7pt]
\prod_{i}{(x^{2^i}-1)}^{j_i}, & \nu_2(k)\ge n-a',~{\rm or}~k=0 .
\end{cases}
\end{equation}

{\em Case 2}: $r=-1+2^{n-a}v$ with $v$ being odd and $\nu_2(\ell)=0$ (i.e. $a'=a$).
By Lemma~\ref{AC^k}(ii) (applied to $A'C^k$),
$${\rm char}_{(A^\ell C^k)|_{V^*}}(x)
=\begin{cases}{\Phi_{2^{a+1}}(x)}^{2^{n-a-1}}, & \mbox{$k$ is odd};\\[7pt]
 {(x^{2^{a}}-1)}^{2^{n-a-1}}, & \mbox{$k$ is even}. \end{cases}
$$
Since $A^\ell C^k$ is similar to a permutation matrix,
by Lemma \ref{permutable} and its corollary,
\begin{equation}\label{r=3 V^2}
{\rm char}_{A^\ell C^k|_{V^2}}(x)=\begin{cases}
{(x^{2^{a}}-1)}^{2^{n-a-1}}, & \mbox{$k$ is odd};\\[7pt]
\prod_{i}{(x^{2^i}-1)}^{j_i}, & \mbox{$k$ is even}.
\end{cases}
\end{equation}

In both cases, by Lemma \ref{permutable} again,
we can conclude that $A^\ell C^k|_{V^2}$ is similar to a permutation matrix.
\qed

\medskip
Now we can get enough information for $A|_{\cal E}$ when
$r{\not\equiv}-1\!\pmod{2^n}$. We proceed in two cases:
${\rm ord}_{{\Bbb Z}_{2^n}^*}(r)=2$, or
 ${\rm ord}_{{\Bbb Z}_{2^n}^*}(r)>2$.

\begin{lemma} \label{a=1}
Let $A,C\in{\rm GL}_{2^n}({\Bbb C})$.
Assume that $C$ is a maximal cycle,
$A^{-1}CA=C^r$ where $-1\ne r\in{\Bbb Z}_{2^n}^*$
with ${\rm ord}_{{\Bbb Z}_{2^{n}}^*}(r)=2$, and $A^2=I$.
Assume further that
$\langle A,C\rangle$ is a permutation-like matrix group.
\begin{itemize}
\item[\rm(i)]
If $r=1+2^{n-1}$, then $A|_{V^2}=I$.
\item[\rm(ii)]
If $r=-1+2^{n-1}$, then
${\rm char}_{A|_{V^2}}(x)=(x-1)^2(x^2-1)^{2^{n-1}-1}$.
\end{itemize}
In the both cases there is a choice of the basis ${\cal E}$ in Eqn \eqref{basis}
such that $A|_{\cal E}$ is a permutation matrix.
\end{lemma}

\pf By Lemma \ref{V^2}, 
$\big\langle A|_{V^2},C|_{V^2}\big\rangle$
is a permutation-like matrix group of dimension $2^{n-1}$
with $C|_{V^2}$ being a maximal cycle.

(i).~
Obviously,  $r=1+2^{n-1}$ acts by multiplication on $2{\Bbb Z}_{2^n}$ trivially.
So, $A|_{{\cal E}^2}$ is a diagonal matrix, hence
$A|_{{\cal E}^2}$ commutes with $C|_{{\cal E}^2}$ which
is a maximal cycle of dimension $2^{n-1}$.
By Lemma \ref{self-centralized},
$A|_{{\cal E}^2}\in\big\langle C|_{{\cal E}^2}\big\rangle$.
But $(A|_{{\cal E}^2})^2=I$ and ${\rm ord}(C|_{{\cal E}^2})=2^{n-1}$.
Thus $A|_{{\cal E}^2}=I_{2^{n-1}\times 2^{n-1}}$
or $= (C|_{{\cal E}^2})^{2^{n-2}}$.
Suppose that $A|_{{\cal E}^2}= (C|_{{\cal E}^2})^{2^{n-2}}$,
i.e.,
$(A|_{{\cal E}^2})(C|_{{\cal E}^2})^{- 2^{n-2}}=I_{2^{n-1}\times 2^{n-1}}$.
Then
$$
{\rm char}_{(A|_{{\cal E}^2})(C|_{{\cal E}^2})^{- 2^{n-2}}}(x)
=(x-1)^{2^{n-1}}.
$$
However, $\nu_2(-2^{n-2})=n-2$, and by Eqn~\eqref{r=5 V^2},
$$
{\rm char}_{(A|_{{\cal E}^2})(C|_{{\cal E}^2})^{- 2^{n-2}}}(x)
={(x^2-1)}^{2^{n-2}}.
$$
We reached an impossible equality
``$(x-1)^{2^{n-1}}={(x^2-1)}^{2^{n-2}}$''.
So it is impossible that $A|_{{\cal E}^2}= (C|_{{\cal E}^2})^{2^{n-2}}$,
and it has to be the case that
 $A|_{{\cal E}^2}=I_{2^{n-1}\times 2^{n-1}}$.
By Lemma \ref{V^*}, we already have a basis ${\cal E}^*$ of $V^*$ such that
$A|_{{\cal E}^*}$ is a permutation matrix. Thus
the basis ${\cal E}$ in Eqn \eqref{basis} can be obtained
such that $A|_{\cal E}$ is a permutation matrix.

(ii).~ Note that $r=-1+2^{n-1}\equiv -1~({\rm mod}~2^{n-1})$
(cf. Eqn \eqref{eq 2Z}), so
$$(A|_{V^2})^{-1}(C|_{V^2})(A|_{V^2})=(C|_{V^2})^{-1}.$$
From Lemma \ref{V^2}(ii) and its Eqn~\eqref{r=3 V^2},
we have seen that
$\big\langle A|_{{\cal E}^2},C|_{{\cal E}^2}\big\rangle$
is a permutation-like matrix group of dimension $2^{n-1}$
with $C|_{{\cal E}^2}$ being a maximal cycle, and
$$
{\rm char}_{A C^k|_{V^2}}(x)=\begin{cases}
{(x^{2}-1)}^{2^{n-2}}, & \mbox{$k$ is odd};\\
\prod_{i}{(x^{2^i}-1)}^{j_i}, & \mbox{$k$ is even}.
\end{cases}
$$
Comparing it with Lemma~\ref{BC^k}, we see that
$(A|_{{V}^2})e_{2^{n-1}}=e_{2^{n-1}}$, and,
with a suitable choice of the basis ${\cal E}^2$,
$A|_{{\cal E}^2}$ is a permutation matrix.
\qed

\begin{lemma}\label{a ge 2}
Assume thar $A,C\in{\rm GL}_{2^n}({\Bbb C})$,
 $C$ is a maximal cycle and
$A^{-1}CA=C^r$ where $r\in{\Bbb Z}_{2^n}^*$ with
${\rm ord}_{{\Bbb Z}_{2^n}^*}(r)=2^a$.
If $a\ge 2$ and $A^{2^a}=I$,
then there is a choice of the basis ${\cal E}$ in Eqn \eqref{basis} such that
$A|_{\cal E}$ is a permutation matrix.
\end{lemma}

\pf We prove it by induction on $n$.
By Remark \ref{Z^*},
$$r=\pm 1+2^{n-a}v,~~ \mbox{with $v$ being odd and $2\le a\le n-2$}.$$
From Lemma \ref{V^*}, we have a basis ${\cal E}^*$ of $V^*$ in
Eqn \eqref{V2V*} such that $A|_{{\cal E}^*}$ is a permutation matrix.
By Lemma~\ref{V^2}, $\big\langle A|_{V^2},C|_{V^2}\big\rangle$
is a permutation-like matrix group of dimension $2^{n-1}$
with $C|_{V^2}$ being a maximal cycle. Let ${A'\!=\!A^{2^{a-1}}}$.
Then $A'^2=I$ and $A'^{-1}CA'=C^{r^{2^{a-1}}}$.
Since $2^{a-1}$ is even, we see that
${r^{2^{a-1}}\!\in\!\langle 5\rangle}$,
hence $r^{2^{a-1}}$ is the unique element of order $2$
in $\langle 5\rangle$, i.e.,
$r^{2^{a-1}}\!\equiv 1+2^{n-1}~({\rm mod}~2^n)$;
see Remark \ref{Z^*}.
By Lemma \ref{a=1}(i), $A'|_{V^2}=I$; i.e.,
$$(A|_{V^2})^{2^{a-1}}=A^{2^{a-1}}|_{V^2}=I.$$
Note that ${\rm ord}_{{\Bbb Z}_{2^{n-1}}^*}(r)=2^{a-1}$.
By induction, we have a basis ${\cal E}^2$ of $V^2$ in Eqn~\eqref{V2V*}
such that $A|_{{\cal E}^2}$ is a permutation matrix.
Then ${\cal E}={\cal E}^*\cup{\cal E}^2$ is a basis of ${\Bbb C}^{2^n}$
in Eqn~\eqref{basis} such that $A|_{{\cal E}}$ is a permutation matrix.
\qed

\begin{theorem}\label{cyclic}
Assume that ${\cal G}$ is a $2^n$-dimensional permutation-like matrix group
which contains a maximal cycle $C$ such that $\langle C\rangle$
is normal in ${\cal G}$ and the quotient group ${\cal G}/\langle C\rangle$ is cyclic.
Then ${\cal G}$ is a permutation matrix group.
\end{theorem}

\pf By the assumption we can assume that
${\cal G}=\langle A,C\rangle$ and
$A^{-1}CA=C^r$ for an $r\in{\Bbb Z}_{2^n}^*$
with ${\rm ord}_{{\Bbb Z}_{2^n}^*}(r)=2^a$ where $0\le a\le n-2$.
By Lemma \ref{self-centralized},
$\langle C\rangle$ is self-centralized in ${\cal G}$.
By Lemma \ref{split} and Lemma \ref{BC^k},
we can further assume that $A^{2^a}=I$ and $1\le a\le n-2$.
By Remark \ref{Z^*}, we prove the theorem in three cases.

{\em Case 1}:~ $r=-1$, hence $a=1$ and $A^2=I$.
By Lemma \ref{BC^k}, $Ae_{2^{n-1}}=\pm e_{2^{n-1}}$.
If $Ae_{2^{n-1}}=-e_{2^{n-1}}$, then we replace $A$ by $AC$.
So we can assume that $Ae_{2^{n-1}}=e_{2^{n-1}}$.
By Lemma \ref{BC^k} (i), there is a basis ${\cal E}$ in
Eqn \eqref{basis} such that $A|_{\cal E}$ is a permutation matrix.
By Lemma \ref{SC}, ${\cal G}$ is a permutation matrix group.

{\em Case 2}:~ $r=\pm 1+2^{n-1}$ (hence $a=1$).
By Lemma \ref{split}(ii), we can assume that $A^2=I$.
By Lemma \ref{a=1},
there is a basis ${\cal E}$ in Eqn \eqref{basis}
such that $A|_{\cal E}$ is a permutation matrix.
By Lemma \ref{SC}, ${\cal G}$ is a permutation matrix group.

{\em Case 3}:~ $r=\pm 1+2^{n-a}v$ with $v$ being odd and $2\le a\le n-2$.
By Lemma \ref{split}(ii), we can assume that $A^{2^a}=I$.
By Lemma \ref{a ge 2}, there is a basis~${\cal E}$ in
Eqn~\eqref{basis} such that $A|_{\cal E}$ is a permutation matrix.
By Lemma \ref{SC}, we are done.
\qed

\section{The quotient ${\cal G}/\langle C\rangle$ is non-cyclic}
\label{non-cyclic case}

In this section we consider any $2^n$-dimensional
permutation-like matrix group ${\cal G}$ with a maximal cycle
$C$ such that $\langle C\rangle$ is normal in ${\cal G}$
and the quotient group ${\cal G}/\langle C\rangle$ is not cyclic,
and prove that ${\cal G}$ is in fact a permutation matrix group.

It is easy to extend Lemma \ref{SC} as follows.

\begin{lemma}
\label{SC 2}
Let $A,B,C\in{\rm GL}_{2^n}({\Bbb C})$ with 
$C$ being a maximal cycle and the subgroup $\langle C\rangle$ 
being normal in ${\cal G}=\langle A,B, C\rangle$.
If there is a choice of the basis ${\cal E}$ in Eqn \eqref{basis}
such that both $A|_{\cal E}$ and $B|_{\cal E}$
are permutation matrices,
then the matrix group ${\cal G}$ is a permutation matrix group.
\end{lemma}

\pf Let $f=\sum_{e\in{\cal E}}e$ and
${\cal F}=\{f,~Cf,~,C^2f,~\cdots,~C^{2^n-1}f\}$.
Because the matrix to transform the basis ${\cal E}$ 
to the set ${\cal F}$ is a 
Vandemond matrix formed by $1,\lambda,\cdots,\lambda^{2^n-1}$
which is non-degenerate, ${\cal F}$ is a basis of ${\Bbb C}^{2^n}$.
With respect to the basis ${\cal F}$, $C|_{\cal F}$ is 
obviously a permutation matrix.
Since $A$ permutes the elements of ${\cal E}$, we have
$Af=f$. For any $C^jf\in{\cal F}$,
$$
A(C^jf)=(AC^jA^{-1})(Af)=C^{r^{-1}j}f\in{\cal F}.
$$
Thus $A|_{\cal F}$ is a permutation matrix. 
Similarly, $B|_{\cal F}$ is a permutation matrix.
So, with respect to the basis ${\cal F}$ 
every element of ${\cal G}$ is a permutation matrix.
\qed

\begin{lemma}\label{r and -1}
Let $r=1+2^{n-a}\in{\Bbb Z}_{2^n}^*$ where $1\le a\le n-2$.
Then there are $j_1,\cdots,j_m\in{\Bbb Z}_{2^n}\backslash\{0,2^{n-1}\}$ such that
all the $\mu_r$-orbits on ${\Bbb Z}_{2^n}\backslash\{0,2^{n-1}\}$
are written as follows (note that $b_t=1$, $1\le t\le m$, is allowed):
$$
\begin{array}{ll}
\Gamma_1=\{j_1,j_1r,\cdots,j_1r^{2^{b_1}-1}\},
    &\mu_{-1}\Gamma_1=\{-j_1,-j_1r,\cdots,-j_1r^{2^{b_1}-1}\},\\[7pt]
\Gamma_2=\{j_2,j_2r,\cdots,j_2r^{2^{b_2}-1}\},
   &\mu_{-1}\Gamma_2=\{-j_2, -j_2r,\cdots, -j_2r^{2^{b_2}-1}\},\\[7pt]
\cdots~~\cdots & \cdots~~\cdots \\[7pt]
 \Gamma_m=\{j_m,j_mr,\cdots,j_mr^{2^{b_m}-1}\},
 & \mu_{-1}\Gamma_m=\{-j_m, -j_mr,\cdots, -j_mr^{2^{b_m}-1}\}.
\end{array}
$$
\end{lemma}

\pf For $k\in{\Bbb Z}_{2^n}$,
$rk\equiv -k~({\rm mod}~2^n)$ if and only if
$(r+1)k\equiv 0~({\rm mod}~2^n)$. However, $\nu_2(r+1)=1$
since $r+1=2+2^{n-a}$ and $n-a\ge 2$. We see that
$rk\equiv -k~({\rm mod}~2^n)$ if and only if
$k\equiv 0~ {\rm or}~2^{n-1}~({\rm mod}~2^n)$.
So, for any $\mu_r$-orbit $\Gamma=\{j,jr,\cdots,jr^{2^{b}-1}\}$
on ${\Bbb Z}_{2^n}\backslash\{0,2^{n-1}\}$, 
the set $\mu_{-1}\Gamma=\{-j, -jr,\cdots, -jr^{2^{b}-1}\}$ is
another $\mu_r$-orbit on ${\Bbb Z}_{2^n}\backslash\{0,2^{n-1}\}$,  and
$$\mu_{-1}(\mu_{-1}\Gamma)
 =\mu_{-1}(\{-j, -jr,\cdots, -jr^{2^{b}-1}\})=\Gamma.$$

Starting from any $\mu_r$-orbit
$\Gamma_1=\{j_1, j_1r,\cdots,j_1r^{2^{b_1}-1}\}$
on ${\Bbb Z}_{2^n}\backslash\{0,2^{n-1}\}$,
we get a pair $\Gamma_1, \mu_{-1}\Gamma_1$ of
$\mu_r$-orbits on ${\Bbb Z}_{2^n}\backslash\{0,2^{n-1}\}$.
Taking any $\mu_r$-orbit
$\Gamma_2=\{j_2, j_2r,\cdots,j_2r^{2^{b_2}-1}\}$  on
${\Bbb Z}_{2^n}\backslash
  \big(\{0,2^{n-1}\}\cup\Gamma\cup\mu_{-1}\Gamma\big)$,
we get another pair $\Gamma_2,\mu_{-1}\Gamma_2$ on
${\Bbb Z}_{2^n}\backslash\{0,2^{n-1}\}$.
Go on in this way till to ${\Bbb Z}_{2^n}\backslash\{0,2^{n-1}\}$ exhausted,
we obtain all the $\mu_r$-orbits on ${\Bbb Z}_{2^n}\backslash\{0,2^{n-1}\}$
as listed pairwise  in the lemma. 
\qed

\smallskip
\begin{theorem}\label{non-cyclic}
Assume that ${\cal G}$ is a $2^n$-dimensional permutation-like matrix group
which contains a maximal cycle $C$ such that $\langle C\rangle$
is normal and the quotient group ${\cal G}/\langle C\rangle$ is not cyclic.
Then ${\cal G}$ is a permutation matrix group.
\end{theorem}

\pf Since $\langle C\rangle$ is self-centralized in ${\cal G}$
(see Lemma \ref{self-centralized}), ${\cal G}/\langle C\rangle$ is isomorphic
to a subgroup ${\cal H}$ of the multiplicative group ${\Bbb Z}_{2^n}^*$.
By Lemma \ref{sub Z^*} (ii), 
${\cal H}=\langle r\rangle\times\langle -1\rangle$
where $r=1+2^{n-a}$ with $1\le a\le n-2$.
Then there are $A,B\in{\cal G}$ such that
$${\cal G}=\langle A,B,C\rangle,~~~
 B^{-1}CB=C^{-1},~~~ A^{-1}CA=C^r.
$$ 
By Lemma~\ref{split}(i) and Lemma \ref{BC^k},
we further assume that $A^{2^a}\!=\!I$ and $B^2\!=\!I$. 

First we show that {\em there is a suitable choice of $B$ such that $AB=BA$}.
Since the quotient group 
${\cal G}/\langle C\rangle\cong {\cal H}$ is commutative,
in ${\cal G}/\langle C\rangle$ the image of
 $B^{-1}AB$ coincides with the image of $A$ itself.
So there is a $C^k\in\langle C\rangle$ such that $B^{-1}AB=AC^k$.
Note that ${\rm ord}(B^{-1}AB)={\rm ord}(A)=2^a$, by Eqn \eqref{(AC)^j},
$$
I=(B^{-1}AB)^{2^a}=(AC^k)^{2^a}=A^{2^a}C^{k(r^{2^a-1}+\cdots+r+1)}
=C^{k(r^{2^a-1}+\cdots+r+1)}.
$$
By Lemma \ref{1+r+}(ii),
$\nu_2(r^{2^a-1}+\cdots+r+1)=a$. So
the above equality implies that $\nu_2(k)\ge n-a$.
We can write $k=2^{n-a}h$. Noting that $r-1=2^{n-a}$, we obtain:
\begin{align*}
 (BC^h)^{-1}A(BC^h)
 &=C^{-h}B^{-1}ABC^h=C^{-h}AC^kC^h=AC^{-rh}C^kC^h\\
 &=AC^{k+h-rh}=AC^{2^{n-a}h-2^{n-a}h} = A.
\end{align*}
Replacing $B$ with $BC^h$, we still have ${\cal G}=\langle A,B,C\rangle$,
$B^{-1}CB=C^{-1}$ and $B^2=I$ 
(since $(BC^h)^2=I$, see Lemma \ref{split}(ii.1)),
and we have a further condition that $AB=BA$.

By Lemma \ref{a=1}(i)  (for $a=1$) and Lemma \ref{a ge 2}  (for $a>1$), 
there exists a basis ${\cal E}$
in Eqn \eqref{basis} such that $A|_{\cal E}$ is a permutation matrix.
Specifically, $Ae_0=e_0$ for any $e_0\in{\rm E}(\lambda^0)$ 
and $Ae_{2^{n-1}}=e_{2^{n-1}}$  for any 
$e_{2^{n-1}}\in{\rm E}(\lambda^{2^{n-1}})$, see Remark~\ref{e_0}.
Using the notation in Lemma \ref{r and -1} and
the method in Eqn \eqref{1-monomial cycle},
we construct the basis ${\cal E}$ precisely by choosing 
$e_j$ for $j\in{\Bbb Z}_{2^n}\backslash\{0,2^{n-1}\}$ as follows.

Let $\Gamma_t$ and $\mu_{-1}\Gamma_t$, for $t=1,\cdots,m$,
be the $\mu_r$-orbits in Lemma \ref{r and -1}.
We take a non-zero vector $e_{j_t}\in{\rm E}(\lambda^{j_t})$. 
By Eqn~\eqref{1-monomial cycle},  
${\cal E}_t=\{e_{j_t}, Ae_{j_t}, \cdots, A^{2^{b_t}-1}e_{j_t}\}$
is a basis of the subspace
$V_t=\bigoplus_{i=0}^{2^{b_t}-1}{\rm E}(\lambda^{r^i j_t})$
such that $A|_{{\cal E}_t}$ is a permutation matrix.
Since $Be_{j_t}\in{\rm E}(\lambda^{-j_t})$ (see Eqn \eqref{permute E}),
we can take $e_{-j_t}=Be_{j_t}$.
Then ${\cal E}'_t=\{e_{-j_t}, Ae_{-j_t}, \cdots, A^{2^{b_t}-1}e_{-j_t}\}$
is a basis of the subspace
$V'_t=\bigoplus_{i=0}^{2^{b_t}-1}{\rm E}(\lambda^{-r^i j_t})$
such that $A|_{{\cal E}'_t}$ is a permutation matrix.
Let ${\cal E}=\{e_0,e_{2^{n-1}}\}\bigcup_{t=1}^m({\cal E}_t\cup{\cal E}'_t)$.
Then ${\cal E}$ is a basis in Eqn~\eqref{basis} such that
$A|_{\cal E}$ is a permutation matrix.

On the other hand, $Be_{-j_t}=B^2e_{j_t}=e_{j_t}$.
i.e., $B$ interchanges $e_{j_t}$ and $e_{-j_t}$.
Since $B$ commutes with $A$,  we get that
$B(Ae_{j_t})=A(Be_{j_t})=Ae_{-j_t}$.
It is the same that $BAe_{-j_t}=Ae_{j_t}$.
So, $B$ interchanges $Ae_{j_t}$ and $Ae_{-j_t}$.
Similarly, $B$ interchanges $A^i e_{j_t}$ and $A^i e_{-j_t}$ for
$i=0,1,\cdots,2^{b_{j_t}-1}$.
In a word, $B|_{{\cal E}_t\cup{\cal E}'_t}$ is a permutation matrix.
Of course, $Be_0=e_0$, see Remark \ref{e_0}.  
To complete the proof,
it is enough to show that ${Be_{2^{n-1}}=e_{2^{n-1}}}$;
because: it implies that  $B|_{\cal E}$ is a permutation matrix,
and by Lemma \ref{SC 2}, ${\cal G}$ is a permutation matrix group.

From an easy computation:
$$(AB)^{-1}C(AB)=B^{-1}(A^{-1}CA)B
 =B^{-1}C^{1+2^{n-a}}B=C^{-1-2^{n-a}},$$
we have $(AB)^{-1}C(AB)=C^{-1+2^{n-a}v}$ with $v=-1$ being odd.
Obviously, $\langle AB, C\rangle$ is a permutation-like matrix group
and $(AB)^{2^a}=A^{2^a}B^{2^a}=I$ (since $AB=BA$).
By Lemma~\ref{a=1}(ii) and Lemma~\ref{a ge 2},
there is a basis ${\cal E}'$ in Eqn~\eqref{basis} such that
$(AB)|_{{\cal E}'}$ is a permutation matrix. In particular,
$(AB)e_{2^{n-1}}=e_{2^{n-1}}$, see Remark \ref{e_0}. 
Note that $Ae_{2^{n-1}}=e_{2^{n-1}}$ implies that
$A^{-1}e_{2^{n-1}}=e_{2^{n-1}}$.
Then $Be_{2^{n-1}}=A^{-1}(ABe_{2^{n-1}})
  =A^{-1}e_{2^{n-1}}=e_{2^{n-1}}$.
We are done.
\qed

\section*{Acknowledgments}

The authors are supported by NSFC through the grant number 11271005.


\begin{thebibliography}{99}

\bibitem{AB}J. L. Alperin, R. B. Bell, Groups and Representations,
GTM 162, Springer-Verlag, New York,  1997.

\bibitem{CDFL}
Bocong Chen, Hai Q. Dinh, Yun Fan, and San Ling,
Polyadic constacyclic codes, 
IEEE Trans. Inform. Theory 61(2015),  4895-4904.

\bibitem{C05} G. Cigler, Groups of matrices with prescribed spectrum,
Doctoral dissertation, 2005,
http://matknjiz.si/doktotati/2005/10921-83.pdf

\bibitem{C07} G. Cigler, Permutation-like matrix groups, Linear Algebra and its
Applications 422(2007), 486-505.

\bibitem{DF} Guodong Deng, Yun Fan, Permutation-like
matrix group with a maximal cycle of prime square length,
Linear Algebra and its Applications 450(2014), 44-55.

\bibitem{DF15} Guodong Deng, Yun Fan,
Permutation-like matrix groups with a maximal cycle of power of odd prime length,
Linear Algebra and its Applications 480(2015), 1-10.


\bibitem{R} D.J.S. Robinson, A Course in Groups Theory,
GTM 80, Springer-Verlag, New York,  1980.

\end{thebibliography}
\end{document}